\documentclass{llncs}
\usepackage{latexsym}
\usepackage{amssymb}
\usepackage{amsfonts}
\usepackage{amsmath}
\usepackage{indentfirst}
\usepackage{graphicx}
\usepackage{epsfig}
\usepackage{slashbox,colortbl}

\newcommand{\cvd}{\hfill $\blacksquare$\bigskip}

\date{}

\author{Lapo Cioni\and Luca Ferrari\inst{}\thanks{Partially supported by INdAM-GNCS 2017 project:
``Codici di stringhe e matrici non sovrapponibili".}}

\institute{Dipartimento di Matematica e Informatica ``U.
Dini"\\
University of Firenze, Firenze, Italy\\
\tt{lapo.cioni@stud.unifi.it, luca.ferrari@unifi.it}}

\title{Enumerative Results on the\\
Schr\"oder Pattern Poset}

\begin{document}

\maketitle

\begin{abstract}
The set of Schr\"oder words (\emph{Schr\"oder language}) is endowed with a natural partial order,
which can be conveniently described by interpreting Schr\"oder words as lattice paths.
The resulting poset is called the \emph{Schr\"oder pattern poset}.
We find closed formulas for the number of Schr\"oder words covering/covered by a given Schr\"oder word
in terms of classical parameters of the associated Schr\"oder path.
We also enumerate several classes of \emph{Schr\"oder avoiding words} (with respect to the length),
i.e. sets of Schr\"oder words which do not contain a given Schr\"oder word.
\end{abstract}

\medskip

%

\medskip

\section{Introduction}

In the literature several definitions of patterns in words can be found.
In the present article we consider a notion of pattern which is rather natural when words are interpreted as \emph{lattice paths},
by using each letter of the alphabet of the word to encode a possible \emph{step}.
The notion of pattern in a lattice paths investigated here has been introduced in \cite{BBFGPW,BFPW}, where it has been studied in the case of Dyck paths.
Aim of the present work is to find some analogous enumerative results in the case of Schr\"oder paths.
In order to make this paper self-contained, we will now briefly recall the main definitions and notations concerning patterns in paths,
and we introduce the basic notions concerning the Schr\"oder pattern poset.

\bigskip

\begin{figure}[!h]
\centering
\includegraphics[scale=0.45]{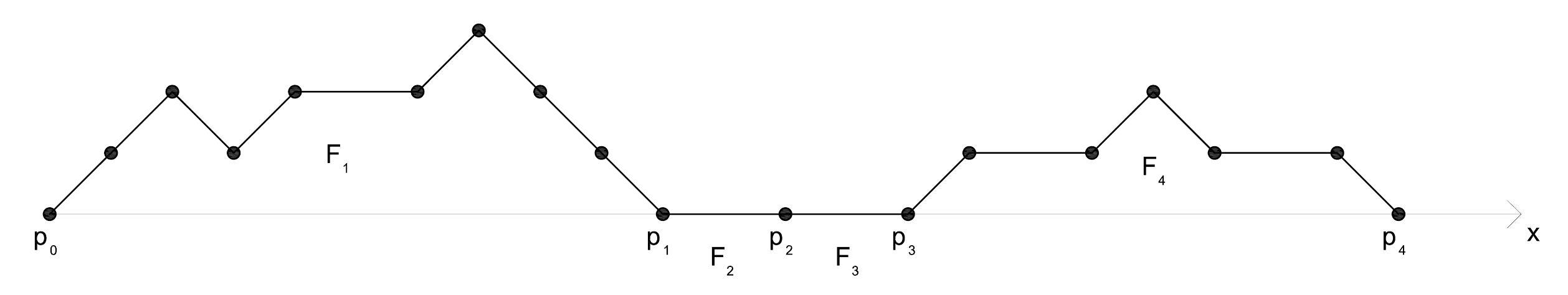}\\
\caption{The Schr\"oder word $UUDUH_2UDDDH_2 H_2 UH_2 UDH_2 D$ represented as a Schr\"oder path.}
\label{Schroder}
\end{figure}

For our purposes, a \emph{lattice path} is a path in the discrete plane starting at the origin of a fixed Cartesian coordinate system,
ending somewhere on the $x$-axis, never going below the $x$-axis and using only a prescribed set of steps $\Gamma$.
We will refer to such paths as \emph{$\Gamma$-paths}.
As a word, a $\Gamma$-path can be represented by the sequence of the letters encoding the sequence of its steps.
In view of this, in the following we will often use the terms ``path" and ``word" referred to the same object.
Classical examples of lattice paths are Dyck, Motzkin and Schr\"oder paths, which are obtained by taking $\Gamma$ to be the set of steps
$\{U,D\}$, $\{U,D,H\}$ and $\{U,D,H_2 \}$, respectively (see Figure \ref{Schroder}).
Here letters represents the steps $U(p)=(1,1)$, $D(own)=(1,-1)$, $H(orizontal)=(1,0)$ and $H_2(\textnormal{\emph{orizontal of length 2}})=(2,0)$, respectively.

Given a $\Gamma$-path $P$, its \emph{length} is given by the final abscissa of $P$.
Also important is the \emph{word length} of $P$, which is the length of the word associated with $P$.
For instance, the Schr\"oder path in Figure \ref{Schroder} has length 22 and word length 17.
Notice that the length of a Schr\"oder path is necessarily even;
for this reason it is sometimes more meaningful to refer to the \emph{semilength} of a Schr\"oder path.

Given two $\Gamma$-paths $P$ and $Q$, we declare $P\leq Q$ whenever $P$ occurs as a (not necessarily contiguous) subword of $Q$.
In this case, we say that $P$ is a \emph{pattern} of $Q$.
So, for instance, the Schr\"oder path $UH_2 UDDH_2 UH_2 H_2 D$ is a pattern of the Schr\"oder path in Figure \ref{Schroder}.
When $P$ is not a pattern of $Q$ we will also say that $Q$ \emph{avoids} $P$.
The set $\mathcal{P}_{\Gamma}$ of all $\Gamma$-paths endowed with the above binary relation is clearly a poset.

\bigskip

In the case of Schr\"oder paths, the resulting poset will be denoted $\mathcal{S}$.
It is immediate to see that $\mathcal{S}$ has a minimum (the empty path), does not have maximum and is locally finite (i.e. all intervals are finite).
Moreover, $\mathcal{S}$ is a ranked poset, and the rank of a Schr\"oder path is given by its semilength.
An important fact concerning $\mathcal{S}$ is that it is a \emph{partial well order},
i.e. it contains neither an infinite properly decreasing sequence nor an infinite antichain
(this is actually a consequence of a well known theorem by Higman \cite{H}).
Notice that this is not the case in another famous pattern poset, the permutation pattern poset, where infinite antichains do exist (see \cite{SB}).

The present paper is devoted to the investigation of some structural and enumerative properties of the Schr\"oder pattern poset.
Specifically, in Section \ref{cover} we study the covering relation of $\mathcal{S}$
and in Section \ref{enumeration} we enumerate some classes of Schr\"oder paths avoiding a single pattern.

We would like to remark that, even when $\Gamma$-paths are interpreted as words over a suitable alphabet,
other kinds of \emph{patterns} can be defined which are equally natural and interesting.
Just to mention one of the most natural, one can require an occurrence of a pattern to be constituted by \emph{consecutive letters}.
This originates what is sometimes called the \emph{factor order}, which has been studied for instance in \cite{B}
(in the unrestricted case of all words on a given alphabet).
Many papers, such as \cite{BGPP}, investigate properties and applications of this more restrictive notion of pattern,
also extending it to the case in which the pattern is not required to be a $\Gamma$-path itself.

\section{The Covering Relation in the Schr\"oder Pattern Poset}\label{cover}

In the Schr\"oder pattern poset $\mathcal{S}$, following the usual notation for the covering relation,
we write $P\prec Q$ ($Q$ \emph{covers} $P$) to indicate that $P\leq Q$ and the rank of $P$ is one less than the rank of $Q$ (i.e., $\text{rank}(P)=\text{rank}(Q)-1$).
The results contained in the present section concern the enumeration of Schr\"oder paths covered by and covering a given Dyck path $Q$.
We need some notation before stating them.

In a Schr\"oder path $Q$, let $k+1$ be the number of points of $Q$ (having integer coordinates) lying on the $x$-axis (call such points $p_0 =(0,0), p_1, \ldots ,p_k$).
Then $Q$ can be factorized (in a unique way) into $k$ Schr\"oder \emph{factors} $F_1 ,\ldots ,F_k$, each $F_i$ starting at $p_{i-1}$ and ending at $p_i$.
Denote with $f_i$ and $h_i$ the number of $U$ and $H_2$ steps of factor $F_i$, respectively. Notice that $f_i$ also equals the number of $D$ steps of the same factor.
Let $a_i$ (resp., $d_i$) be the number of \emph{ascents} (resp. \emph{descents}) in $F_i$, where an ascent (resp. descent) is a maximal consecutive run of $U$ (resp., $D$) steps.
Moreover, we denote with $p(Q)$ and $v(Q)$ the number of occurrences in $Q$ of a consecutive factor $UDU$ and $DUD$, respectively.
Finally, we denote with $h(Q)$ the total number of \emph{flats} of $Q$, a flat being a maximal sequence of consecutive $H_2$ steps.
The path depicted in Figure \ref{Schroder} has 4 factors, and we have $f_1 =4, f_2 =f_3 =0, f_4 =2$, $h_1 =h_2 =h_3 =1, h_4 =2$, $a_1 =3, a_2 =a_3 =0, a_4 =2$,
$d_1 =2, d_2 =d_3 =0, d_4 =2$, $p(Q)=1$, $v(Q)=0$ and $h(Q)=4$.

\begin{proposition}
If $Q=F_1 F_2 \cdots F_k$ is a Schr\"oder path with $k$ factors, with
$F_i$ having $a_i$ ascents and $d_i$ descents, then the number of Schr\"oder paths covered
by $Q$ is given by
\begin{equation}\label{covered}
\sum_{1\leq i\leq j\leq k} d_i a_j -p(Q)-v(Q)+h(Q)\quad .
\end{equation}
\end{proposition}

\emph{Proof.}\quad There are two (mutually exclusive) ways to obtain a Schr\"oder path covered by $Q$, namely:
\begin{enumerate}
\item by removing a $H_2$ step, or
\item by removing a $U$ step and a $D$ step.
\end{enumerate}

We examine the two cases separately.

\begin{enumerate}
\item It is immediate to observe that one obtains the same path by removing any of the steps belonging to the same flat,
whereas removing a step from different flats gives rise to different paths.
Therefore, the number of distinct Schr\"oder paths obtained by $Q$ by removing an $H_2$ step is $h(Q)$.
\item We wish to prove that there are $\sum_{1\leq i\leq j\leq k} d_i a_j -p(Q)-v(Q)$ ways to remove a $U$ steps and a $D$ steps from $Q$
and to obtain another Schr\"oder path. We will proceed by induction on the number $k$ of factors of $Q$.
If $k=1$, then necessarily $Q$ starts with a $U$ step and ends with a $D$ step (otherwise $Q=H_2$, which has no $U$ and $D$ steps).
Observe that, in this case, we can remove any of the $U$ steps and any of the $D$ steps and the resulting path is still a Schr\"oder path.
Removing steps from the same ascent (and from the same descent) returns the same path,
so we have $a_1$ possible choices to remove a $U$ step and $d_1$ possible choices to remove a $D$ step.
However, there are some special cases in which, though removing from different ascents or descents, we obtain the same path.
Specifically, if we have a consecutive string $UDU$ in $Q$, then removing from $Q$ the $UD$ of such a string returns the same path as removing the $DU$,
in spite of the fact that the two $U$ steps belong to different ascents.
In a similar way, the presence of a factor $DUD$ in $Q$ gives the possibility of getting the same Schr\"oder paths by removing steps from different descents.
To avoid overcount, we thus have to subtract the number of consecutive strings $UDU$ and $DUD$ of $Q$, thus obtaining a total of $d_1 a_1-p(Q)-v(Q)$ paths.
Now suppose that $k>1$. There are three distinct cases to analyze.

\begin{itemize}
\item If we remove both the $U$ and the $D$ steps from the prefix $F_1 \cdots F_{k-1}$ of $Q$ consisting of the first $k-1$ factors
(which is a Schr\"oder path in itself, of course), by induction we have $\sum_{1\leq i\leq j\leq k-1} d_i a_j -p(F_1 \cdots F_{k-1})-v(F_1 \cdots F_{k-1})$
distinct choices.
\item If we remove both the $U$ and the $D$ steps from the last factor $F_k$, using the same argument as the case $k=1$ we get
$d_k a_k -p(F_k )-v(F_k )$ distinct paths.
\item Finally, suppose we choose to remove the $D$ step from the prefix $F_1 \cdots F_{k-1}$ and the $U$ step from the last factor $F_k$
(notice that we are not allowed to do the opposite, otherwise the resulting path would not be Schr\"oder).
In this case we have $a_k$ possible choices for $U$ and $\sum_{i=1}^{k-1}d_i$ possible choices for $D$.
Once again, however, there are some paths that are overcounted, occurring when $F_1 \cdots F_{k-1}$ and $F_k$ share a consecutive $UDU$ or a consecutive $DUD$.
A quick look shows that this overcount is corrected by subtracting the number of such shared occurrences of consecutive $UDU$ and $DUD$.
\end{itemize}

The sum of the above three cases gives the required expression.

\end{enumerate}

Finally, summing up the quantities obtained in 1. and 2., we obtain precisely formula (\ref{covered}).\cvd

\emph{Remark.}\quad If $Q$ is a Dyck paths, then $h(Q)=0$, and formula (\ref{covered}) reduces to the analogous formula for Dyck paths obtained in \cite{BFPW,BBFGPW},
since a Schr\"oder path covered by a Dyck path is necessarily a Dyck path.

\begin{proposition}
Let $P=F_1 \cdots F_k$ be a Schr\"oder path having $k$ factors. Denote with $f_i$ the number of $U$ steps in the factor $F_i$
(this is also the number of $D$ steps in $F_i$) and with $h_i$ the number of $H_2$ steps in $F_i$. Moreover, let $\ell$ be the word length of $P$.
Then the number of Schr\"oder paths covering $P$ is given by
\begin{equation}\label{covering}
2+\ell +\sum_{(i,j)\atop 1\leq i\leq j\leq k}(f_i +h_i )(f_j +h_j ).
\end{equation}
\end{proposition}

\emph{Proof.}\quad We have two options to obtain a Schr\"oder path covering $P$:
\begin{enumerate}
\item either we add a $H_2$ step, or
\item we add a $U$ step and a $D$ step.
\end{enumerate}

As in the previous proposition, we examine the two cases separately.

\begin{enumerate}
\item Adding a new $H_2$ step in any point of a flat of $P$ returns the same path.
Hence, in order to obtain distinct paths, we can add a $H_2$ step either before a $U$ step, or before a $D$ step, or at the end of $P$.
Thus we have a total of
$$
2\sum_{i=1}^{k}f_i +1
$$
paths covering $P$ in this case.
\item We start by observing that adding a $U$ step in any point of an ascent returns the same path, and the same holds for $D$ steps (with ascents replaced by descents).
Suppose to add a new $U$ step to $P$ first. In order to obtain distinct paths, we can add $U$ either before a $D$ step, or before a $H_2$ step, or at the end of $P$.

If a $U$ step is added before a $D$ step in $F_i$, we observe that we cannot add the new $D$ step in a factor $F_j$, with $j<i$,
otherwise the path would fall below the $x$-axis.
With this constraint in mind, we are now allowed to add the new $D$ step either before a $U$ or before a $H_2$ or at the end of the path.
However, in the first of the three previous cases, we cannot of course add the new $D$ step before the first allowed $U$ (i.e., at the beginning of the factor);
moreover, adding the new $D$ right before the new $U$ step just added would produce a substring $DUD$,
which can be obtained also by first adding the $U$ step in the following available position of $P$ and then adding the $D$ step immediately after it.
Thus, in this case, the number of paths covering $P$ is obtained by considering the number of possible choices for $U$ to be added in $F_i$, which is $f_i$,
and the number of possible choices for $D$, which is $\sum_{j\geq i}(f_j +h_j )$, and so it is
$$
\sum_{(i,j)\atop 1\leq i\leq j\leq k}f_i (f_j +h_j ).
$$

If a $U$ step is added before a $H_2$ step in $F_i$, as in the previous case, we cannot add the new $D$ step in a factor $F_j$, with $j<i$.
We can now add the new $D$ step either before a $U$ (except for the first $U$ of $F_i$, of course), or before a $H_2$ or at the end of $P$.
So, in this case, the number of paths covering $P$ is given by
$$
\sum_{i=1}^{k}h_i \cdot \left( 1+\sum_{j=i}^{k}(f_j +h_j )\right) .
$$

Finally, if we add the new $U$ step at the end of $P$, then the new $D$ step must necessarily be added after it.

\end{enumerate}

Summing up, we therefore obtain the following expression for the total number of paths covering $P$:
\begin{align}\label{last}
&2\sum_{i=1}^{k}f_i +1+\sum_{(i,j)\atop 1\leq i\leq j\leq k}f_i (f_j +h_j )+\sum_{i=1}^{k}h_i \cdot \left( 1+\sum_{j=i}^{k}(f_j +h_j )\right) +1 \nonumber\\
=& 2+\sum_{i=1}^{k}(2f_i +h_i )+\sum_{(i,j)\atop 1\leq i\leq j\leq k}f_i (f_j +h_j )+\sum_{(i,j)\atop 1\leq i\leq j\leq k}h_i (f_j +h_j )\nonumber \\
=& 2+\ell +\sum_{(i,j)\atop 1\leq i\leq j\leq k}(f_i +h_i )(f_j +h_j ),
\end{align}
which is formula (\ref{covering}).\cvd

\emph{Remark 1.}\quad Notice that $f_i +h_i$ is the semilength of the factor $F_i$.
Denoting it with $\varphi_i$, formula (\ref{covering}) can be equivalently written as
$$
2+\ell +\sum_{(i,j)\atop 1\leq i\leq j\leq k}\varphi_i \varphi_j .
$$

\emph{Remark 2.}\quad If $P$ is a Dyck path, then, in the first expression in the chain of equalities (\ref{last}), the summand $2\sum_{i=1}^{k}f_i +1$
gives the number of non-Dyck paths covering $P$ (i.e., those having one $H_2$ step), and the remaining summands give the number of Dyck paths covering $P$.
Also in this case, recalling that $h_i =0$ for all $i$, we recover the analogous formula for Dyck paths obtained in \cite{BFPW,BBFGPW}.

\section{Enumerative Results on Pattern Avoiding Schr\"oder Paths}\label{enumeration}

Main goal of the present section is to enumerate several classes of Schr\"oder paths avoiding a given pattern.
For any Schr\"oder path $P$, denote with $\mathcal{S}_n (P)$ the set of Schr\"oder paths of semilength $n$ avoiding $P$,
and let $s_n (P)=|\mathcal{S}_n (P)|$ be its cardinality. It is completely trivial to observe that
\begin{itemize}
\item $s_n (\emptyset )=0$;
\item $s_n (H_2 )=C_n$, where $C_n =\frac{1}{n+1}{2n\choose n}$ is the $n$-th Catalan number (sequence A000108 of \cite{S}), counting the number of Dyck paths of semilength $n$;
\item $s_n (UD)=1$ when $n>0$.
\end{itemize}

Starting from patterns of semilength 2, we get some interesting enumerative results.
In the next subsections we define several classes of Schr\"oder paths avoiding a single pattern,
each of which suitably generalizes the case of a pattern of semilength 2.
In all cases, after having described a general enumeration formula, we illustrate it in the specific case of the relevant pattern of semilength 2.

Before delving into computations we state an important lemma, which in several cases reduces the enumeration of pattern avoiding Schr\"oder paths
to the case of pattern avoiding Dyck paths. In this lemma, as well as in several subsequent proofs, we will deal with the \emph{multiset coefficient}
$\left( {n\choose k}\right)$, counting the number of multisets of cardinality $k$ of a set of cardinality $n$. As it is well known,
the multiset coefficients can be expressed in terms of the binomial coefficients, namely $\left( {n\choose k}\right) ={n+k-1\choose k}$.

\begin{lemma}\label{dyck}
Given a Dyck path $P$, denote with $d_n (P)$ the number of Dyck paths of semilength $n$ avoiding $P$.
Then
\begin{equation}\label{avoid_dyck}
s_n (P)=\sum_{h=0}^{n}{n+h\choose n-h}d_h (P).
\end{equation}
\end{lemma}

\emph{Proof.}\quad Let $Q$ be a Schr\"oder path.
Clearly $Q$ avoids $P$ if and only if the Dyck path $\tilde{Q}$ obtained from $Q$ by deleting horizontal steps avoids $P$.
Therefore, the set of Schr\"oder paths of semilength $n$ with $n-h$ horizontal steps avoiding $P$ is obtained by
taking in all possible way a Dyck path of semilength $h$ avoiding $P$ and then adding to it $n-h$ horizontal steps in all possible ways.
Observe that, in a Dyck path of semilength $h$, one has $2h+1$ possible sites where to insert a horizontal step,
and any number of horizontal steps can be inserted into the same site.
Thus, if $n-h$ horizontal steps have to be inserted, it is necessary to select a multiset of cardinality $n-h$ from the set of the possible $2h+1$ sites.
This can be done in ${(2h+1)+(n-h)-1\choose n-h}={n+h\choose n-h}$ ways, as it is well known.
Since $h$ can be chosen arbitrarily in the set $\{ 0,1,2,\ldots ,n\}$, the total number of Schr\"oder paths of semilength $n$ avoiding $P$ is given by
formula (\ref{avoid_dyck}), as desired.\cvd

\subsection{The Pattern $(UD)^k$}

Since $(UD)^k$ is a Dyck path, this case can be seen as an immediate consequence of Lemma \ref{dyck} together with the results of \cite{BFPW,BBFGPW}.

\begin{proposition}
For $i,j\geq 1$, let $N_{i,j}=\frac{1}{i}{i\choose j}{i\choose j-1}$ be the Narayana numbers (sequence A001263 of \cite{S}).
Extend such a sequence by setting $N_{0,0}=1$ and $N_{i,0}=N_{0,j}=0$, for all $i,j>0$. Then
\begin{equation}\label{(UD)^k}
s_n ((UD)^k )=\sum_{h=0}^{n}\sum_{j=0}^{k-1}{n+h\choose n-h}N_{h,j}.
\end{equation}
\end{proposition}

The case $k=2$ gives rise to an interesting situation.
In fact, for the pattern $UDUD$, recalling that $N_{0,0}=1$, $N_{i,0}=0$ and $N_{i,1}=1$ for all $i>0$, formula (\ref{(UD)^k}) gives
\begin{align}\label{fibonacci}
s_n (UDUD)&=\sum_{h=0}^{n}{n+h\choose n-h}N_{h,0}+\sum_{h=0}^{n}{n+h\choose n-h}N_{h,1}\nonumber \\
&=1+\sum_{h=1}^{n}{n+h\choose n-h}=\sum_{h=0}^{n}{2n-h\choose h}.
\end{align}

Since it is well known that Fibonacci numbers $(F_n)_n$ (sequence A000045 in \cite{S}) can be expressed in terms of binomial coefficients
as\footnote{Notice that the sum contains only a finite number of nonzero terms, since the binomial coefficients vanish when $k>\lfloor n/2 \rfloor$}
$F_{n+1} =\sum_{k\geq 0}{n-k\choose k}$,
we get $s_n (UDUD)=F_{2n+1}$, i.e. Schr\"oder paths avoiding $UDUD$ are counted by Fibonacci numbers having odd index (sequence A122367 of \cite{S}).

\bigskip

\emph{Remark.}\quad Notice that, for a Schr\"oder path, avoiding the (Schr\"oder) path $UDUD$ is equivalent to avoiding the (non-Schr\"oder) path $DU$.
As suggested by one of the referees, a simple combinatorial argument to count Schr\"oder words of semilength $n$ avoiding the subword $DU$ is the following:
if the words contains $k$ $H_2$ steps, then it can be constructed by taking the word $U^{n-k}D^{n-k}$ and inserting $k$ $H_2$ steps.
Taking into account all possibilities, and summing over $k$, gives precisely formula (\ref{fibonacci}).


\subsection{The Pattern $U^k D^k$}\label{casoUkDk}

This case is similar to the previous one, in that it can be easily inferred from Lemma \ref{dyck}, since the generic pattern of the class is a Dyck path.
Thus, applying the above mentioned lemma and using results of \cite{BFPW,BBFGPW}, we obtain the following result.

\begin{proposition}
For all $k\geq 0$, we have
\begin{align}\label{U^kD^k}
s_n (U^k D^k )&=\sum_{h=0}^{k-1}{n+h\choose n-h}C_h +{n+k\choose n-k}(C_k -1)\nonumber \\
&+\sum_{h=k+1}^{\min \{ 2k-1,n\} }\sum_{j\geq 1}{n+h\choose n-h}b_{k-j,h-k+j}^{2},
\end{align}
where the $C_n$'s are the Catalan numbers and the $b_{i,j}$'s are the ballot numbers (sequence A009766 of \cite{S}).
\end{proposition}

Setting $k=2$ in formula (\ref{U^kD^k}), for Schr\"oder paths of semilength $n\geq 3$ avoiding $UUDD$ we obtain the following polynomial of degree 4
(sequence A027927 of \cite{S}):
$$
s_n (UUDD)=1+{n+1\choose n-1}+{n+2\choose n-2}=1+\frac{n(n+1)(n^2 +n+10)}{24}.
$$


\subsection{The Pattern $H_2 ^k$}

In this case the generic pattern of this class is not a Dyck path. However, we are able to give a direct argument to count Schr\"oder paths avoiding $H_2 ^k$.

\begin{proposition}
For all $n,k\geq 0$, we have
\begin{equation}\label{H_2^k}
s_n (H_2 ^k )=\sum_{i=0}^{k-1}{2n-i\choose i}C_{n-i} .
\end{equation}
\end{proposition}

\emph{Proof.}\quad We observe that a Schr\"oder path avoids the pattern $H_2 ^k$ if and only if it has at most $k-1$ $H_2$ steps.
Thus, the set of all Schr\"oder paths of semilength $n$ avoiding $H_2 ^k$ can be obtained by taking the set of all Dyck paths of semilength $n-i$
and inserting in all possible ways $i$ $H_2$ steps, for $i$ running from $0$ to $k-1$.
Since in a Dyck path of semilength $n-i$ there are precisely $2n-2i+1$ points in which inserting the horizontal steps, and we have to insert $i$ horizontal steps
(possibly inserting more than one $H_2$ step in the same place), we get
\begin{align*}
s_n (H_2 ^k )&=\sum_{i=0}^{k-1}\left( {2n-2i+1\choose i} \right) C_{n-i}=\sum_{i=0}^{k-1}{2n-2i+1+i-1\choose i}C_{n-i}\\
&=\sum_{i=0}^{k-1}{2n-i\choose i}C_{n-i},
\end{align*}
as desired.\cvd

When $k=2$ we obtain the following special case:
$$
s_n (H_2 H_2)=C_n +(2n-1)C_{n-1}=\frac{n+3}{2}C_n =\frac{n+3}{2(n+1)}{2n\choose n},
$$
which is valid for $n\geq 1$. This is sequence A189176 of \cite{S}, whose generating function is
$\frac{1-5x+4x^2 -(1-5x)\sqrt{1-4x}}{2x(1-4x)}$, and can be also obtained as the row sums of a certain Riordan matrix (see \cite{S} for details).

\subsection{The Pattern $UH_2 ^{k-1}D$}

This class of patterns requires a little bit more care, nevertheless we are able to get a rather neat enumeration formula.

\begin{proposition}
For all $n\geq 0,k>1$, we have
\begin{equation}\label{UH_2^{k-1}D}
s_n (UH_2 ^{k-1}D)=1+\sum_{h=1}^{n}\sum_{i=0}^{\min \{ k-2,n-h\}}\left( {2h-1\choose i}\right) (n-h-i+1)C_h .
\end{equation}
\end{proposition}

\emph{Proof.}\quad Let $P$ be a Schr\"oder path of semilength $n$ avoiding $UH_2 ^{k-1}D$.
If $P$ does not contain $U$ steps, then necessarily $P=H_2 ^n$. Otherwise, $P$ can be decomposed into three subpaths, namely:
\begin{itemize}
\item a prefix, consisting of a (possibly empty) sequence of horizontal steps;
\item a path starting with the first $U$ step and ending with the last $D$ step (necessarily not empty);
\item a suffix, consisting of a (possibly empty) sequence of horizontal steps.
\end{itemize}

The central portion of $P$ in the above decomposition has at most $k-2$ horizontal steps.
Thus it can be obtained starting from a Dyck path of semilength $h$, for some $1\leq h\leq n$, and adding $i$ horizontal steps, for some $0\leq i\leq \min \{ k-2,n-h\}$.
Such horizontal steps can be inserted into $2h-1$ possible sites
(here we have to exclude the starting and the ending points, since the subpath is required to start with a $U$ and end with a $D$),
with the possibility of inserting several steps into the same site, as usual.
The resulting path has therefore semilength $h+i$, and has to be completed by adding a suitable number of horizontal steps at the beginning and at the end,
to obtain a Schr\"oder path of semilength $n$: there are $n-h-i+1$ possible ways to do it.
We thus obtain formula (\ref{UH_2^{k-1}D}) for $s_n (UH_2 ^{k-1}D)$, as desired.\cvd

Formula (\ref{UH_2^{k-1}D}) becomes much simpler in the special case $k=2$ of Schr\"oder paths of semilength $n$ avoiding $UH_2 D$. Indeed we obtain:
$$
s_n (UH_2 D)=1+\sum_{h=1}^{n}(n-h+1)C_h =\sum_{h=0}^{n}C_h (n-h)+\sum_{h=0}^{n}C_h -n.
$$

In this case, we can find an interesting expression for the generating function of these coefficients
in terms of the generating function $C(x)=\sum_{n\geq 0}C_n x^n$ of Catalan numbers, which provides an easy way to compute them:
\begin{align*}
\sum_{n\geq 0}s_n (UH_2 D)x^n &=C(x)\cdot \sum_{n\geq 0}nx^n +C(x)\cdot \sum_{n\geq 0}x^n -\sum_{n\geq 0}nx^n \\
&=C(x)\left( \frac{x}{(1-x)^2}+\frac{1}{1-x} \right) -\frac{x}{(1-x)^2}\\
&=\frac{1}{(1-x)^2}(C(x)-x).
\end{align*}

Roughly speaking, the above generating function tells us that $s_n (UH_2 D)$ is given by the partial sums of the partial sums of the sequence of Catalan numbers
where $C_1$ is replaced by 0. The associated number sequence starts with $1,2,5,13,35,99,295,\ldots$ and does not appear in \cite{S}.

\subsection{The Pattern $H_2 ^{k-1}UD$}

The last class of patterns we consider is the most challenging one. It gives rise to an enumeration formula which is certainly less appealing than the previous ones.
Due to space limitation, we will just sketch its proof and simply state the special case corresponding to $k=2$.
Before illustrating our final results, we need to introduce a couple of notations.

\bigskip

We denote with $P_{k,h}$ the number of Dyck prefixes of length $k$ ending at height $h$.
Notice that we can express these coefficients in terms of the ballot numbers $b_{i,j}=\frac{i-j+1}{i+1}{i+j\choose i}$,
counting the number of Dyck prefixes with $i$ up steps and $j$ down steps, as follows:
$$
P_{h,k}=\sum_{i,j\atop {i-j=h\atop i+j=k}}b_{i,j}.
$$

Moreover, we denote with $S_{n,q}$ the number of Schr\"oder paths of semilength $n$ having exactly $q$ horizontal steps.
Since each such path can be uniquely determined by a Dyck path of semilength $n-q$ with $q$ horizontal steps added,
we have a rather easy way to compute $S_{n,q}$ in terms of the Catalan numbers $C_n$:
$$
S_{n,q}=C_{n-q}\left( {2(n-q)+1\choose q} \right) ={2n-q\choose q}C_{n-q}.
$$

\begin{proposition}
For all $n\geq 0,k>1$, we have
\begin{align}\label{H_2^{k-1}UD}
s_n (H_2 ^{k-1}UD)&=\sum_{q=0}^{k-2}S_{n,q}\nonumber \\
&+\sum_{p=0}^{2n-2k+2}\sum_{h=0}^{2n-p-2k+2}P_{p,h}\left( {p+1\choose k-2}\right) \left( {\frac{2n-h-p-2k+4}{2}\choose h}\right) .
\end{align}
\end{proposition}

\emph{Proof.}\quad Let $P$ be a Schr\"oder path of semilength $n$. If $P$ contains less than $k-1$ horizontal steps, then it necessarily avoids $H_2 ^{k-1}UD$.
This gives the term $\sum_{q=0}^{k-2}S_{n,q}$ in the r.h.s. of (\ref{H_2^{k-1}UD}).
If instead $P$ contains at least $k-1$ horizontal steps, then, in order to avoid $H_2^{k-1}UD$, it has to be decomposable into
a Schr\"oder prefix ending at some height $h$ and having exactly $k-2$ horizontal steps, and a suffix starting with a horizontal steps followed exclusively by $H_2$ and $D$ steps
(with exactly $h$ $D$ steps).
The generic Schr\"oder prefix of the required form can be obtained by taking a Dyck prefix of length $p$, for some $0\leq p\leq 2n-2k+2$,
and adding $k-2$ $H_2$ steps in all possible ways.
We thus get a total of $P_{p,h}\cdot \left( {p+1\choose k-2}\right)$ Schr\"oder prefixes of length $p+2k-4$ ending at height $h$ and containing exactly $k-2$ $H_2$ steps.
The generic suffix of the required form contains $h$ $D$ steps and $\frac{2n-h-p-2k+4}{2}$ $H_2$ steps.
Such a suffix can be obtained by inserting the $h$ down steps into the sequence of horizontal steps in all possible ways.
Since the first step has to be a horizontal one, this gives a total of $\left( {\frac{2n-h-p-2k+4}{2}\choose h}\right)$ allowed suffixes.
Putting together all the contributions, we get the desired expression for $s_n (H_2 ^{k-1}UD)$.\cvd

Specializing to $k=2$ we obtain
$$
s_n (H_2 UD)=C_n +\sum_{p=0}^{2n-2}\sum_{h=0}^{2n-p-2}P_{p,h}\left( {\frac{2n-h-p}{2}\choose h} \right) .
$$

\section{Suggestions for Further Work}

It would be very interesting to investigate in more detail the structural properties of the Schr\"oder pattern poset.
A typical question in this context concerns the computation of the M\"obius function, which is still open even in the Dyck pattern poset.
Another (partially related) issue is the enumeration of (saturated) chains.
More generally, can we say anything about the order structure of intervals (for instance, is it possible to determine when they are lattices?)?

Concerning the enumeration of pattern avoiding classes, the next step would be to count classes of Schr\"oder words simultaneously avoiding two or more patterns.

Finally, it would be nice to have information on the asymptotic behavior of integer sequences counting pattern avoiding Schr\"oder words.
In the Dyck case, all the sequences which count Dyck words avoiding a single pattern $P$ have the same asymptotic behavior
(which is roughly exponential in the length of $P$). This is in contrast, for instance with the permutation pattern poset,
where the asymptotic behavior of a class of pattern avoiding permutations depends on the patterns to be avoided
(this is the ex Stanley-Wilf conjecture, proven by Marcus and Tardos \cite{MT}).
What does it happen in the Schr\"oder pattern poset?

\end{document}